\documentclass[11pt]{article}
\usepackage[cp1251,cp866]{inputenc}
\usepackage[english]{babel}
\usepackage{amsfonts}
\usepackage[fleqn]{amsmath}
\usepackage{amssymb}
\usepackage{amsthm}
\usepackage{graphicx}
\usepackage{caption2}
\usepackage{algorithm}
\usepackage[usenames]{color}
\usepackage{mathrsfs}

\newtheorem{theorem}{\quad Theorem}[section]
\newtheorem{lemma}{\quad Lemma}[section]
\newtheorem{corollary}{\quad Corollary}[section]
\newtheorem{proclaim}{\quad Proposition}[section]

\newtheorem{conj}{\quad Conjecture}[section]

\renewcommand{\subsubsection}{\@startsection{subsubsection}{3}{\parindent}{3.25ex plus 1ex minus 1.2ex}{1.5ex plus .2ex}{\large\bf}}
\renewcommand{\subsection}{\@startsection{subsection}{2}{0pt}{3.25ex plus 1ex minus 1.2ex}{1.5ex plus .2ex}{\large\bf}}

\begin{document}

\begin{center}
{\LARGE Reconstruction of $p$-disconnected graphs}

\vspace{5mm}

\large

{ Pavel Skums}

\vspace{5mm}

Mechanics and Mathematics faculty \\
Belarus State University\\
Minsk, Belarus\\
Email: skumsp@gmail.com

\vspace{5mm}

April 24, 2008

\vspace{5mm}

\small {\bf Abstract}

\end{center}

We prove that Kelly-Ulam conjecture is true for $p$-disconnected
graphs.

\vspace{3mm}


\vspace{3mm}

Subj-class: Combinatorics

\section{Introduction}
Let $G$ be a simple graph. The collection $D(G)=(G_v)_{v\in V(G)}$
of vertex-deleted subgraphs of graph $G$ is called the {\it deck}
of $G.$ The graph $H$ with deck $D(H)=(H_u)_{u\in V(H)}$ is called
the {\it reconstruction} of $G$ if there exists a bijection
$f:V(G)\rightarrow V(H)$ such that $G_v\cong H_{f(v)}.$ In this
case we say that the decks $D(G)$ and $D(H)$ are {\it equal.}  The
graph $G$ is {\it reconstructible} if it is isomorphic to any of
its reconstructions.

The following conjecture, first posed in 1942, is one of the most
famious open problems in graph theory.
\begin{conj}{Kelly-Ulam reconstruction conjecture} \cite{Kelly42},\cite{Ulam60}.
Every graph with at least three vertices is reconstructible.
\end{conj}

It is clear that a graph is reconstructible if and only if its
complement is reconstructible.

The class of graphs is called {\it reconstructible} if all graphs
from this class with at least three vertices are reconstructible.
The known examples of reconstructible classes are disconnected
graphs, complements of disconnected graphs, regular graphs etc.
(see, for example, \cite{Bondy91}, \cite{BH77}). Analogously, a
graph parameter is {\it reconstructible}, if it is the same for
all graphs with equal decks. For example, it is easy to show that
degree sequence of graph is reconstructible (\cite{Bondy91},
\cite{BH77}).

The class of graphs $\mathcal{R}$ is called {\it recognizable} if
for any graph $G\in \mathcal{R}$ all its reconstructions also
belong to $\mathcal{R}$. The class $\mathcal{R}$ is {\it weakly
reconstructible} if for any $G\in \mathcal{R}$ every
reconstruction of $G$ which belongs to $\mathcal{R}$ is isomorphic
to $G$. Clearly $\mathcal{R}$ {\it is reconstructible if and only
if it is recognizable and weakly reconstructible.}

 We write
 $u\sim v$ $(u\not\sim v)$ if vertices $u$ and $v$ are
 adjacent (non-adjacent). For the subsets $U, W\subseteq V(G)$ the
 notation $U\sim W$ means that $u\sim w$ for all vertices $u\in U$
 and $w\in W$, $U\not\sim W$ means that there are no adjacent vertices
 $u\in U$ and $w\in W$. To shorten notation, we
 write $u\sim W$ $(u\not\sim W)$ instead of $\{u\}\sim W$
 $(\{u\}\not\sim W)$.

 A {\it triad} is a triple $T=(G,A,B)$, where $G$ is a graph and $(A,B)$ is an
 ordered partition of $V(G)$ into two disjoint subsets. Isomorphism of two triads $T=(G,A,B)$ and $S=(H,C,D)$
 is an isomorphism of graphs $G$ and $H$ preserving corresponding
 partitions. In this case we say, that the triads $T$ and $S$ are
 {\it isomorphic} ($T\cong S$).

 Let $G$ be a graph,
 $M\subseteq V(G)$. $M$ is called a {\it module}
 of $G$ if $v\sim M$ or $v\not\sim M$ for every vertex
 $v\in V(G)\setminus M$. If $M$ is a module, then $V(G)$ is naturally
 partitioned into three parts:
\begin{equation}\label{parts}
  V(G)=A\cup B\cup M,\ A\sim M,\ B\not\sim M.
\end{equation}
 The partition (\ref{parts}) is {\it associated} with the module
 $M$. In this case we write $G=T\circ F$, where $T =(G[A\cup
 B],A,B)$, $F\cong G[M]$.

  For every graph $G$ the sets $V(G)$, singleton subsets of $V(G)$ and
 $\emptyset$ are modules. The modules $M$ with $1<|M|<|VG|$ are
 called {\it nontrivial modules} or {\it homogeneous sets}.

 A graph $G$ is called {\it 1-decomposable} \cite{TC85b}, if there exists a
 module $M$ (called ({\it 1-module})) of $G$ with associated partition $(A,B,M)$ such that
 $A$ is a clique and $B$ is a stable set. Otherwise $G$ is called {\it 1-indecomposable}.
  The properties and applications of
 1-decomposable graphs are described, for example, in
\cite{MP95},\cite{Tysh00},\cite{BLS99}. One of the most important
for us facts concerning 1-decomposable graphs is the following
result of V. Turin.

\begin{theorem}\cite{Turin87} 1-decomposable graphs are reconstructible
\end{theorem}

A graph $G$ is called $P_4$-{\it connected} (or $p$-{\it
connected}), if for every partition of $V(G)$ into two disjoint
sets $V_1$ and $V_2$ there exists an induced $P_4$ (called {\it
crossing} $P_4$) which contains vertices from both $V_1$ and
$V_2$. Otherwise $G$ is called $P_4$-{\it disconnected} (or
$p$-{\it disconnected}). $P_4$-disconnected graphs were introduced
by B. Jamison and S. Olariu in  \cite{JO95}.  The $p$-{\it
connected component} of $G$ is a maximal induced $p$-connected
subgraph of $G$. It is clear that every disconnected graph is
$p$-disconnected, but inverse inclusion is not true.

A graph is called {\it split} \cite{HS81}, if there exists a
partition of its set of vertices $V(G)=A\cup B$ into a clique and
a stable set. This partition is called a {\it bipartition} and
denoted as $(A,B)$.

In this paper we prove that $P_4$-disconnected graphs are
reconstructible. In particular, it generalizes the results about
reconstructibility of disconnected graphs, complements of
disconnected graphs and 1-decomposable graphs.

Let $A$ be a subset of vertices of $G$ such that $G[A]\cong P_4$.
A {\it partner} of $A$ in G is a vertex $v \in G\setminus A$ such
that $G[A\cup v]$ contains at least two induced $P_4$s. A graph
$G$ is $P_4$-{\it tidy} \cite{GRH97}, if any $P_4$ has at most one
partner. The class of $P_4$-tidy graphs contains well-know classes
of $P_4$-extensible, $P_4$-lite, $P_4$-reducible, $P_4$-sparse,
$P_4$-free graphs (see \cite{GRH97}).

We show that the reconstructibility of $P_4$-disconnected graphs
implies the reconstructibility of $P_4$-tidy graphs. Therefore, in
particular, all listed above classes are also reconstructible.
Note, that the reconstructibility of $P_4$-reducible graphs was
proved by B. Thatte in \cite{Thatte95}.

\section{Reconstruction of $p$-disconnected graphs.}
A $p$-connected graph $S$ is called {\it separable} \cite{JO95},
if there exists a disjoint partition of its vertex set $V(S) =
A\cup B$ such that every crossing $P_4$ has its midpoints in $A$
and its endpoints in $B$. In this case a triad $(S,A,B)$ is called
a {\it separable p-connected triad}.

\begin{lemma}\cite{JO95}\label{uniqpartsep}
Every separable $p$-connected graph induces a unique separable
$p$-connected triad.
\end{lemma}

Let's call a triad $(G,A,B)$ {\it generalized split triad}, if
every connected component of $\overline{G[A]}$ and $G[B]$ is a
module in $G$. For example, if all connected components of
$\overline{G[A]}$ and $G[B]$ consist of one vertex, then $G$ is a
split graph.

\begin{lemma}\cite{JO95}\label{strucsep}
Let $T = (G,A,B)$ be separable $p$-connected triad. Then $T$ is a
generalized split triad. Moreover, the graphs $\overline{G[A]}$,
$G[B]$ are disconnected.
\end{lemma}
Note, that, in particular, separable $p$-connected triad contains
at least four vertices.

A split graph $G$ with bipartition $(A,B)$ is called {\it spider},
if there exists a bijection $f:B\rightarrow A$ such that one of
the following conditions holds:

\begin{itemize}
\item[1)] $N(b) = \{f(b)\}$ for every vertex $b\in B$ ({\it thin
spider});

\item[2)] $N(b) = A\setminus \{f(b)\}$ for every vertex $b\in B$
({\it thick spider}).
\end{itemize}

\begin{theorem}\cite{HS81}\label{splitgraphs}
Let $G$ be a graphs, $V(G)=\{v_1,...,v_n\}$, $deg(v_1)\geq
deg(v_2)\geq \dots\ \geq deg(v_n)$ and let $m=m(G) = max\{i :
deg(v_i)\geq i-1\}$. Then $G$ is split if and only if

\begin{equation}\label{HSineq}
    \sum_{i=1}^m deg(v_i) = m(m-1) + \sum_{i=m+1}^n deg(v_i).
\end{equation}

Moreover, if (\ref{HSineq}) holds, then $A=\{v_1,...,v_m\}$ is a
maximal clique and $B=\{v_{m+1},...,v_n\}$ is a stable set.
\end{theorem}

\begin{lemma}\label{reconspiders}
Spiders are reconstructible.
\end{lemma}

\begin{proof} Since thick spiders are complements of thin spiders, it is
sufficient to prove that thin spiders are reconstructible.

Let $G$ be a graph with $V(G)=\{v_1,...,v_n\}$, $deg(v_1)\geq
deg(v_2)\geq \dots\ \geq deg(v_n)$. Taking into account Theorem
\ref{splitgraphs}, it is evident that $G$  is a thin spider if and
only if (\ref{HSineq}) and the following conditions hold:

\begin{itemize}
\item[1)] $deg(v_i) = m(G)$ for every $i=1,\dots,m(G)$;

\item[2)] $deg(v_i) = 1$ for every $i=m(G)+1,\dots,n$.
\end{itemize}

Since degree sequence of graph is reconstructible, thin spiders
are reconstructible. \end{proof}

A vertex $v$ in a $p$-connected graph $G$ is called {\it
$p$-articulation vertex}, if $G_v$ is $p$-disconnected. If every
vertex of $G$ is a $p$-articulation vertex, then $G$ is called
{\it minimally $p$-connected.}

\begin{theorem}\cite{Bab98,BabOl98}\label{minpconspider}
Graph $G$ is minimally $p$-connected if and only if $G$ is a
spider.
\end{theorem}

\begin{theorem}\cite{Bab98}\label{twoparticul}
A $p$-connected graph which is not minimally $p$-connected
contains at least two vertices which are not $p$-articulation
vertices.
\end{theorem}

The following structure theorem was proved in \cite{JO95}. In our
terms it could be written in the following way:

\begin{theorem}\cite{JO95}.\label{structthm}
For an arbitrary graph $G$ exactly one of the following statements
is true:

\begin{enumerate}
\item [1)] $G$ is disconnected;

\item[2)] $\overline{G}$ is disconnected ($G$ is
antidisconnected);

\item[3)] there is a unique separable component $S$ of $G$ with
corresponding partition $V(S) = A\cup B$ such that $G =
(S,A,B)\circ H$;

\item[4)] $G$ is $p$-connected.
\end{enumerate}
\end{theorem}

For example, all connected and anticonnected $1$-decomposable
graphs satisfy 3).

Let $\mathcal{R}$ be the class of graphs $G$ such that

\begin{enumerate}
\item[a)] $G$ is $p$-disconnected;

\item[b)] $G$ is both connected and anticonnected;

\item[c)] $G$ is 1-indecomposable.

\end{enumerate}

To prove, that $p$-disconnected graphs are reconstructible, by
Theorem \ref{structthm} it is sufficient to prove that class
$\mathcal{R}$ is reconstructible.

\begin{lemma}\label{multgenspgrnotpcon}
Let $T$ be generalized split triad and let $H$ be an arbitrary
graph. Then $G = T\circ H$ is $p$-disconnected.
\end{lemma}
\begin{proof} Let $V(G)=A\cup B \cup C$ such that $(G[A\cup
B],A,B)\cong T$, $G[C]\cong H$ and $G = (G[A\cup B],A,B)\circ
G[C]$. It is easy to see that for the partition

\begin{equation}\label{partitionnotp4}
    (A\cup B, C)
\end{equation}

there is no crossing $P_4$. Indeed, let vertices $x,y,z,t$ induces
crossing $P_4$ for the partition (\ref{partitionnotp4}) with
midpoints $y,z$ and endpoints $x,t$ such that $y\sim x$, $z\sim
t$. The only possibility is $x\in C$, $y\in A$, $z,t\in B$. Then
the vertices $z$ and $t$ belongs to the same connected component
$U$ of $S[B]$. But since $U$ is a homogeneous set and $y\sim z$ we
have $y\sim t$. The contradiction is obtained.\end{proof}

As a corollary we obtain that 1-decomposable graphs are
$p$-disconnected.

\begin{lemma}\label{pdisonecard}
Graph is $p$-disconnected if and only if it is not a spider and at
most one of its cards is $p$-connected.
\end{lemma}
\begin{proof} Assume, that $G$ is $p$-disconnected graph. By Theorem
\ref{minpconspider} $G$ is not a spider. Let's show that at most
one card of $G$ is $p$-connected.

If $G$ ($\overline{G}$) is disconnected, then clearly at most one
card of $G$ is connected (anticonnected), therefore our statement
is true. Let $G =T\circ H$, where $T = (S,A,B)$ is separable
$p$-connected triad. If $|H| > 1$, then all cards of $G$ has the
form $T_v\circ H$ or $T\circ H_v$. Thus by Lemma
\ref{multgenspgrnotpcon} all cards of $G$ are $p$-disconnected. If
$|H| = 1$, then $D(G)=\{T_v\circ H\}\cup\{S\}$. Therefore by Lemma
\ref{multgenspgrnotpcon} there exists the unique $p$-connected
card of $G$, isomorphic to $S$.

Inversely, let $G$ is not a spider and at most one of its card is
$p$-connected. Suppose that $G$ is $p$-connected. Then by Theorem
\ref{twoparticul} there exist at least two $p$-connected cards of
$G$. This is contradiction. \end{proof}

Since spiders are reconstructible, the following corollary is
true.

\begin{corollary}\label{recognnonpconnect}
$p$-disconnected graphs are recognizable
\end{corollary}

Since disconnected graphs, antidisconnected graphs and
1-decomposable graphs are reconstructible, we have

\begin{corollary}\label{recognR}
Class $\mathcal{R}$ is recognizable
\end{corollary}

In the further considerations we will use the following technical
lemma.

\begin{lemma}\label{multgenspgrnotpconcomponent}
Let $G = (G[A\cup B],A,B)\circ G[C]$, where $(G[A\cup B],A,B)$ is
generalized split triad, and let $D$ be $p$-connected component of
$G$. Then $D\subseteq A\cup B$ or $D\subseteq C$.
\end{lemma}

\begin{proof} Suppose that $D\cap (A\cup B)\ne \emptyset$, $D\cap C\ne
\emptyset$. As it was shown in Lemma \ref{multgenspgrnotpcon}, for
the partition $(A\cup B, C)$ there is no crossing $P_4$ in $G$.
Therefore for the partition

\begin{equation}
    (D\cap (A\cup B),D\cap C)
\end{equation}

there is no crossing $P_4$ in the graph $G[D]$. This contradicts
the fact, that $D$ is $p$-connected component of $G$. \end{proof}

\begin{lemma}\label{weakrecpdis}
The class $\mathcal{R}$ is weakly reconstructible.
\end{lemma}

\begin{proof} Let $G^1=T^1\circ H^1$, $G^2=T^2\circ H^2$ be two graphs from
$\mathcal{R}$ with equal decks $D(G^1)$ and $D(G^2)$,
$T^1=(S^1,A^1,B^1)$, $T^2=(S^2,A^2,B^2)$ are separable
$p$-connected triads from the definition of the class
$\mathcal{R}$. By Theorem \ref{structthm} $G^1\cong G^2$ if and
only if $T^1\cong T^2$ and $H^1\cong H^2$.

Let $|H^1| = 1$. Then $D(G^1) = \{T^1_v\circ H^1\}\cup \{S^1\}$.
It is evident, that all vertex-deleted triads $T^1_v$, $T^2_u$ are
generalized split triads. Therefore by Lemma
\ref{multgenspgrnotpcon} there exists a unique $p$-connected card
of $G^1$, and this card is isomorphic to $S^1$.

If $|H^2| > 1$, then $D(G^2) = \{T^2_v\circ H^2\}\cup \{T^2\circ
H^2_u\}$ and hence by Lemma \ref{multgenspgrnotpcon} all cards of
$G^2$ are $p$-disconnected.

Therefore $|H^2|=1$ and there exists a unique $p$-connected card
of $G^2$, isomorphic to $S^2$. Thus we have $S^1\cong S^2$. By
Lemma \ref{uniqpartsep} $T^1\cong T^2$ and consequently $G^1\cong
G^1$.

Let further $|H^1|\geq 2,$ $|H^2|\geq 2$. Assume that $V(G^i) =
A^i\cup B^i\cup C^i$, where $(G[A^i\cup B^i],A^i,B^i)\cong T^i$,
$G^i[C^i]\cong H^i$ and $G^i\cong (G[A^i\cup B^i],A^i,B^i)\circ
G^i[C^i]$, $i=1,2$.

Then

\begin{equation}\label{decksGi}
    D(G^i) = D_{T^i}\cup D_{H^i},
\end{equation}

where

\begin{equation}\label{TiHi}
D_{T^i}=\{T^i\circ H^i_v : v\in C^i\}, D_{H^i} = \{T^i_v\circ H^i
: v\in A^i\cup B^i\}, i=1,2.
\end{equation}

By Lemma \ref{multgenspgrnotpcon} all cards from $D(G^i)$,
$i=1,2$, are $p$-disconnected. Clearly all cards from $D_{T^i}$,
$i=1,2$ are both connected and anticonnected $p$-disconnected
graphs (since so is $G^i$).

\begin{proclaim}\label{propforweakpdis}
Let $G^1_v=T^1\circ H^1_v\in D_{T^1}$, $G^2_u=T^2_u\circ H^2\in
D_{H^2}$ and $G^1_v\cong G^2_u$. Then $|T^1| < |T^2|.$
\end{proclaim}

\begin{proof}  Put $C^1_v = C^1\setminus \{v\}$, $A^2_u = A^2\setminus
\{u\}$.

Let $\varphi : V(G^1)\setminus\{v\}\rightarrow V(G^2)\setminus
\{u\}$ be isomorphism of graphs $G^1_v$ and $G^2_u$. If
$\varphi(A^1\cup B^1)\subseteq C^2$, then $\varphi(C^1_v)\supseteq
A^2_u\cup B^2$. But then, for example, $B^2\sim C^2\cap
\varphi(A^1)$, that is impossible.

Therefore by Lemma \ref{multgenspgrnotpconcomponent} it is true,
that $\varphi(A^1\cup B^1)\subseteq (A^2_u\cup B^2)$. Thus
$|T^1|\leq |T^2_u| < |T^2|.$ \end{proof}

Now let's show that there exist $v\in V(G^1)$ and $u\in V(G^2)$
such that

\begin{equation}\label{uvsimHpdis}
    G_v^1\in D_{T^1},\,\,\, G_u^2\in D_{T^2},\,\,\,
G_v^1\cong G_u^2.
\end{equation}

Suppose the contrary. Then there exist $G_{v_1}^1\in D_{T^1}$,
$G_{v_2}^1\in D_{H^1}$, $G_{u_1}^2\in D_{T^2}$, $G_{u_2}^2\in
D_{H^2}$ such that

\begin{equation}\label{eqT1pdis}
T^1\circ H_{v_1}^1 = G_{v_1}^1 \cong G_{u_2}^2 = T_{u_2}^2\circ
H^2,
\end{equation}
\begin{equation}\label{eqT2pdis}
T^2\circ H_{u_1}^2 = G_{u_1}^2 \cong G_{v_2}^1 = T_{v_2}^1\circ
H^1.
\end{equation}

By Proposition \ref{propforweakpdis}

\begin{equation}\label{ineq1pdis}
    |T^1|  < |T^2|,
\end{equation}

and

\begin{equation}\label{ineq2pdis}
    |T^2|  < |T^1|.
\end{equation}

The contradiction is obtained.

So, consider $v\in V(G^1)$, $u\in V(G^2)$ such that
(\ref{uvsimHpdis}) holds.  We have

\begin{displaymath}
T^1\circ H^1_v\cong T^2\circ H^2_u.
\end{displaymath}

By Theorem \ref{structthm} it is true that

\begin{equation}\label{triadseqpdis}
T^1\cong T^2.
\end{equation}

In particular, if $G^1_v\in D_{H^1}$ and $G^1_v\cong G^2_u$ then
$G^2_u\in D_{H^2}$. Indeed, if there exist the cards from
$D_{T^1}$ and $D_{H^2}$ such that (\ref{eqT1pdis}) holds, then the
inequality (\ref{ineq1pdis}) is true. This contradicts
(\ref{triadseqpdis}).

It remains to prove that $H^1\cong H^2$.

Since $G$ is 1 indecomposable, we have that $S^1$ is not a split
graph. Thus there exists a connected component $X$ of
$\overline{G^1[A^1]}$ or $G^1[B^1]$ such that $|X|> 2$. Therefore
it is easy to see, that for any $v\in X$ $T^1_v$ is separable
$p$-connected triad and the card $G^1_v$ is both connected and
anticonnected $p$-disconnected graph.

Let $v\in X$ and $T^1_v\circ H^1 = G^1_v\cong G^2_u = T^2_u\circ
H^2$ and let $\psi$ be isomorphism of graphs $G^1_v$ and $G^2_u$.
By the same reasoning, as in the proof of Proposition
\ref{propforweakpdis} we have $\psi((A^1\cup B^1)\setminus
\{v\})\subseteq (A^2\cup B^2)\setminus \{u\}$. Since $|T^1| =
|T^2|$, it is true that $\psi((A^1\cup B^1)\setminus \{v\}=
(A^2\cup B^2)\setminus \{u\}$. Therefore $\psi(C^1) = C^2$ and
thus $H^1\cong H^2$. \end{proof}

So, Corollary \ref{recognR} and Lemma \ref{weakrecpdis} imply

\begin{theorem}\label{reconpdis}
$p$-disconnected graphs are reconstructible.
\end{theorem}

A {\it quasi-starfish} (resp. {\it quasi-urchin}) \cite{GRH97} is
a graph obtained from a thick spider (resp. thin spider) by
replacing at most one vertex by a $K_2$ or a $O_2$.

\begin{theorem}\cite{GRH97} A graph $G$ is $P_4$-tidy if and only if every
p-component of $G$ is isomorphic to either a $P_5$ or a
$\overline{P_5}$ or a $C_5$ or a quasi-starfish or a quasi-urchin.
\end{theorem}

\begin{corollary} $P_4$-tidy graphs are reconstructible.
\end{corollary}

\begin{proof} Let $G$ be $P_4$-tidy graph. If $G$ is $p$-disconnected, then
by Theorem \ref{reconpdis} $G$ is reconstructible. Suppose that
$G$ is $p$-connected. Then $G$ is isomorphic to either a $P_5$ or
a $\overline{P_5}$ or a $C_5$ or a quasi-starfish or a
quasi-urchin. Clearly $P_5$, $\overline{P_5}$, $C_5$ are
reconstructible. Moreover, quasi-starfishes are complements of
quasi-urchins and by Lemma \ref{reconspiders} spiders are
reconstructible. Thus it is sufficient to consider the case, then
$G$ is obtained from a thin spider $H$ with bipartition $(A,B)$
and with at least 6 vertices by replacing a vertex $v\in V(H)$ by
$K_2$ or $O_2$. Consider the following cases:

1) $v\in A$ is replaced by $K_2$. In \cite{Tysh00} the complete
description of the structure of 1-indecomposable split unigraphs
is presented. From that description one can see that $G$ is a
split unigraph and thus $G$ is reconstructible.

2) $v\in B$ is replaced by $O_2$. By the same description from
\cite{Tysh00} $G$ is a split unigraph and therefore $G$ is
reconstructible.

3) $v\in B$ is replaced by $K_2$. It is easy to see that a graph
$F$ is isomorphic to $G$ if and only if $|V(F)|=|V(G)|$, there
exist exactly two vertices $x,y\in V(F)$ with $deg(x)=deg(y)= 2$
and $F_x\cong F_y$ is a thin spider. Therefore it is evident, that
$G$ is reconstructible.

4) $v\in A$ is replaced by $O_2$. Then it is also easy to see that
a graph $F$ is isomorphic to $G$ if and only if
$|V(F)|=|V(G)|=2k+1$, $k\geq 3$ and there exist two vertices
$x,y\in V(F)$ such that $deg(x)=deg(y) = k$ and cards $F_x, F_y$
are thin spiders. Thus in this case $G$ is also reconstructible.
\end{proof}


\begin{thebibliography}{00}\label{ispist}

\bibitem{BabOl99}
Babel, L. and S. Olariu. On the $p$-connectedness of graphs - a
survey // Discrete Appl. Math.~--- 1999.~--- Vol. 95.~--- P.
11---33.

\bibitem{Bab98}
Babel, L. Tree-like $P_4$-connected graphs, Discrete Math.~---
1998.~--- Vol. 191.~--- P. 13---23.

\bibitem{BabOl98}
Babel, L. and S. Olariu. On the isomorphism of graphs with few
$P_4$s
// Discrete Appl. Math.~--- 1998.~--- Vol. 84.~--- P. 1---3.

 \bibitem{Bondy91}
Bondy, A.  A graph reconstructor's manual, Surveys in
Combinatorics 1991, London Math. Soc. Lecture Notes Series,
Cambridge University Press, Cambridge, 1991, pp. 221-252.

\bibitem{BH77}
Bondy, A and R.L. Hemminger. Graph reconstruction - a survey, J.
Graph Theory 1 (1977), 227-268.

\bibitem{BLS99}
 Brandst\"adt A., Le V.B. and Spinrad J.
 Graph classes: a survey.~---
 Philadelphia: SIAM monographs on discrete mathematics and
 applications, 1999.


\bibitem{FH77}
 F\"oldes S. and Hammer P.L.
 Split graphs //
 Proceedings of the 8-th South-East Conf. of Combinatorics,
 Graph Theory and Computing.~--- 1977.~--- Vol. 19.~--- P. 311---315.

 \bibitem{GRH97} Giakoumakis V., Roussel F. and Thuillier H.
On $P_4$--tidy graphs // Discrete Math. and Theor. Comp. Sci. 1
1997. P. 17--41.

 \bibitem{HS81}
 Hammer P.L. and Simeone B.
 The splittence of a graph //
 Combinatorica.~--- 1981.~--- Vol. 3, No 1.~--- P. 275---284.

\bibitem{JO95}
 Jamison B., and S. Olariu, p-components and the homogeneous
decomposition of graphs, SIAM Journal of Discrete Math. 8 (1995),
448 - 463.

\bibitem{Kelly42}
Kelly, P.J. On isometric transformations, PhD thesis, University
of Wisconsin (1942).

\bibitem{MP95}
Mahadev, N.V. and U.N. Peled, Threshold graphs and related topics,
Annals of Discrete Math. 56 (1995).

\bibitem{Thatte95} Thatte, B. Some results on the reconstruction
problem. $p$-claw-free, chordal and $P_4$-reducible graphs // J.
of Graph Theory, Vol. 19, No. 4, 549-561 (1995)

\bibitem{Turin87}
Turin V. Reconstruction of decomposable graphs // Vestsi AN BSSR.
1987. N. 3. P. 16-20.

  \bibitem{Tysh00}
 Tyshkevich R.
 Decomposition of graphical sequences and unigraphs //
 Discrete Math.~--- 2000.~--- Vol. 220.~--- P. 201---238.

  \bibitem{TC85b}
 Tyshkevich R. and Chernyak A.
 Decomposition of graphs //
 Cybernetics.~--- 1985.~--- Vol. 21.~--- P. 231---242.

\bibitem{Ulam60} Ulam S.M. A collection of mathematical problems, Wiley
(Interscience), New York. 29 (1960). MR 22:10844
\end{thebibliography}
\end{document}